\documentclass[10pt,twoside]{article}
\usepackage{amssymb}
\usepackage{graphicx}
\usepackage{latexsym}
\newtheorem{theorem}{Theorem}[section]

\newtheorem{proposition}{Proposition}[section]

\def\proof{\mbox {\it Proof.~}}
\def\endproof{\mbox {$\Box$}}


\begin{document}

\title{
\begin{flushleft}
{\small{\it prepint}  (2013)}
\end{flushleft}
\vspace{0.5in}
{\bf\Large  Fat handles and phase portraits of Non Singular Morse-Smale flows on $S^{3}$ with unknotted saddle orbits.}}
\author{{\bf\large B. Campos,}\footnote{The authors acknowledge the support of Ministerio de Ciencia y Tecnolog\'{\i}a MTM2011-28636-C02-02
and the Universitat Jaume I P11B2011-30}\hspace{2mm}
{\bf\large P. Vindel}\vspace{1mm}\\
{\it\small Instituto Universitario de Matem\'{a}ticas y Aplicaciones de
Castell\'{o}n, IMAC.}\\ {\it\small  Departamento de Matem\'{a}ticas. Universitat Jaume I},
{\it\small Castell\'on. 12071. Spain}\\
{\it\small e-mail: campos@uji.es, vindel@uji.es}}\vspace{1mm}
\date {October 2013}
\maketitle

\begin{center}
{\bf\small Abstract}
\vspace{3mm}
\hspace{.05in}\parbox{4.5in}
{{\small In this paper we build Non-singular Morse-Smale flows on S$^{3}$ with
unknotted and unlinked saddle orbits by identifying fat round handles along
their boundaries. This way of building the flows enables to get their phase
portraits.

We also show that the presence of heteroclinic trajectories imposes an order
in the round handle decomposition of these flows; this order is total for
NMS flows composed of one repulsive, one attractive and $n$ unknotted saddle
orbits, for $n\geq 1$.
 }}
\end{center}

\section{\bf Introduction}

The phase portrait of a dynamical system provides a complete description of
the flow. It is not easy to get this picture for 3-dimensional manifolds and
the drawing becomes more difficult when the number of periodic orbits
increases. Our goal is to obtain an easy way for visualizing the flows.

Our study on $S^{3}$ is based on the round handle decomposition of the
manifold for a given flow by D. Asimov \cite{asimov} and by J. Morgan \cite%
{morgan} and the topological characterization of the links of periodic
orbits by M. Wada \cite{wada}.

Every non-singular Morse Smale flow on $S^{3}$ admits a round handle
decomposition whose core circles are the periodic orbits of the flow, and
conversely, a round handle decomposition gives rise to a Non-singular
Morse-Smale flow on a flow manifold.

The characterization of the links of periodic orbits of a non singular
Morse-Smale flow on $S^{3}$ by M. Wada is obtained from Hopf links by
applying six operations; each operation corresponds to a different
attachment of a round 1-handle in the decomposition of the manifold. The
flows characterized by the first three operations have all their saddle
orbits unknotted and unlinked.

In this paper, we consider flows coming only from the first three operations
of Wada. A description of this type of flows can be found in \cite{campos1}
and \cite{campos2}.

To claim our goal we define in section \ref{seccion asas} the basic fat
handles (fat handles with only one saddle orbit) and we prove that these
flows \ can also be obtained from the identification of basic fat handles
along their boundaries (Propositions \ref{nus} and \ref{lnus}). In section %
\ref{sectorder}, we show that the appearance of heteroclinic trajectories
connecting saddle orbits implies noncommutativity of the operations
involved, establishing an order in the round handle decomposition
(Proposition \ref{orden}). We find a total or linear order when the flow
contains only one attractive and one repulsive periodic orbit because of the
existence of transversal intersections of invariant manifolds of the saddle
orbits (Theorem \ref{orden3}).

\section{\bf Previous Results\label{seccionprviousresults}}

A non singular Morse-Smale flow (or NMS for short) is a flow without fixed
points, consisting of a finite number of hyperbolic periodic orbits where
the intersections of stable and unstable manifolds of the saddle orbits are
transversal.

D. Asimov \cite{asimov} and J.W. Morgan \cite{morgan} established a
correspondence between NMS flows and round handle decompositions of the
corresponding manifold. These flows are defined on flow manifolds.

Let $M$ be a compact manifold whose boundary has been partitioned into two unions of components: $\partial M=\partial_{-}M \cup \partial_{+}M$, $\partial_{-}M \cap \partial_{+}M=\emptyset$.
A flow manifold  is a pair $(M, \partial_{-}M)$ satisfying:
\begin{itemize}
  \item $\chi(\partial_{-}M)=\chi(M)$.
  \item $\chi(\partial_{+}M)=\chi(M)$.
  \item There exists a nonsingular vector field on $M$ pointing inwards on $\partial_{-}M$ and outwards on $\partial_{+}M$.
\end{itemize}

For the case of dimension $3$, a pair $\left( M,\partial _{-}M\right) $ of a manifold $M$ and a
compact submanifold $\partial _{-}M$ of $\partial M,$ or by abuse of
notation, a manifold $M$ is called:

A round 0-handle if $\left( M,\partial _{-}M\right) \cong \left(
D^{2}\times S^{1},\emptyset \right).$

A round 1-handle if $\left( M,\partial _{-}M\right) \cong \left(
D^{1}\times D^{1}\times S^{1},D^{1}\times \partial D^{1}\times
S^{1}\right) .$

A round 2-handle if $\left( M,\partial _{-}M\right) \cong \left(
D^{2}\times S^{1},\partial D^{2}\times S^{1}\right).$

In this case, the round handles are diffeomorphic to tori
and correspond to 0-handles when there is a repulsive periodic orbit in the
core, to 2-handles if there is an attractive periodic orbit in the core and
to 1-handles if the orbit in the core is a saddle; $0,$ $1$ and $2$
are the indices of the periodic orbits. A set of indexed periodic orbits is called an indexed link.

\begin{figure}[h!]
\begin{center}
 \includegraphics[scale=0.5]{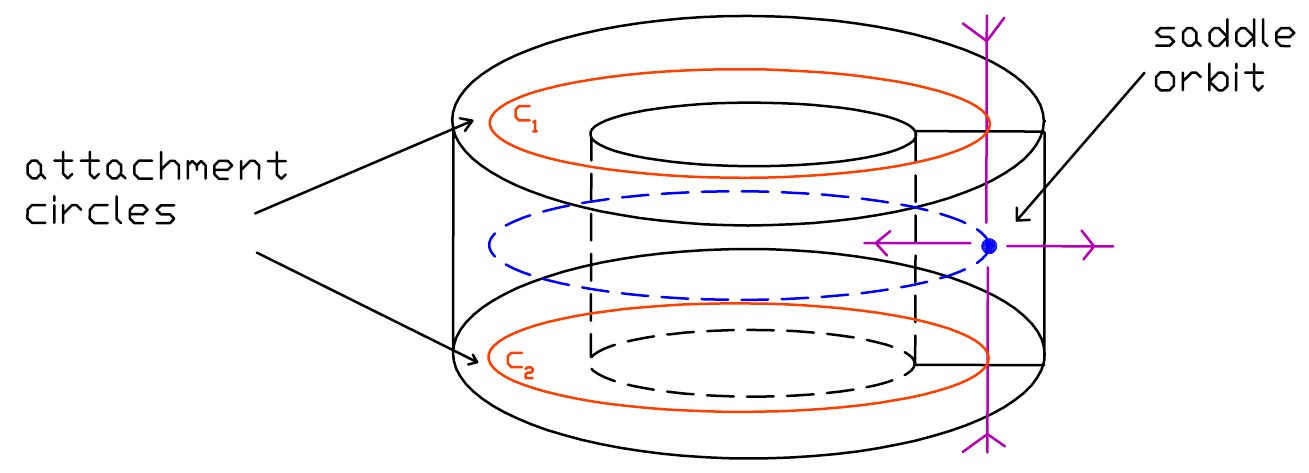}\\
 \caption{Round 1-handle.}\label{1asa}
\end{center}
\end{figure}

\begin{proposition}[Morgan]
Given a flow manifold (X,$\partial _{-}$X) with a NMS flow, then (X,$%
\partial _{-}$X)\ has a round handle decomposition whose core circles are
the closed orbits of the flow.
\end{proposition}

The round handle decomposition for a compact, orientable 3-manifold $M$ was
modified by Morgan:
\begin{equation}
\emptyset =M_{0}\subset M_{1}\subset \cdots \subset M_{i}\subset
M_{i+1}\subset \cdots \subset M_{N}=M  \label{filtra1}
\end{equation}%
where each manifold $M_{i}$, called \emph{fat round handle}, is obtained
from $M_{i-1}$ by attaching a round $1-$handle by means of one or two
attaching circles (see Figure \ref{1asa}).

From the round handle decomposition of the 3-dimensional sphere, \ M. Wada
characterizes the set of periodic orbits of NMS flows in terms of a
generator, the Hopf link with indices $0$ and $2$ attached to the components, and six operations defined from the type of
attachment of the round 1-handles (see \cite{wada}).

\begin{theorem}[M. Wada]
\textquotedblright Every indexed link which consists of all the closed
orbits of a Non-Singular Morse-Smale flow on $S^{3}$  is obtained
from Hopf links by applying the following six operations. Conversely, every
indexed link obtained from Hopf links by applying the operations is the set
of all the closed orbits of some Non-Singular Morse-Smale flow on $S^{3}"$.
\end{theorem}

\emph{OPERATIONS: For given indexed links }$l_{1}$\emph{\ and }$l_{2}$\emph{%
, the six operations are defined as follows. Let }$l_{1}\cdot l_{2}$\emph{\
denote the split sum of the links }$l_{1}$\emph{\ and }$l_{2}$\emph{\ and }$%
N(k,M)$\emph{\ a regular neighborhood of }$k$\emph{\ in }$M$\emph{.}

1) $I(l_{1},l_{2})=l_{1}\cdot l_{2}\cdot u$\emph{, where }$u$\emph{\ is an
unknot with index 1.}

2) $II(l_1,l_2)=l_1\cdot (l_2-k_2)\cdot u$\emph{, where }$k_2$\emph{\ is a
component of }$l_2$\emph{\ of index 0 or 2.}

3) $III(l_{1},l_{2})=(l_{1}-k_{1})\cdot (l_{2}-k_{2})\cdot u$\emph{, where }$%
k_{1}$\emph{\ is a component of }$l_{1}$\emph{\ of index 0 and }$k_{2}$\emph{%
\ is a component of }$l_{2}$\emph{\ of index 2.}

4) $IV(l_1,l_2)=(l_1\#l_2)\cup m$\emph{. The connected sum }$(l_1\#l_2)$%
\emph{\ is obtained by composing a component }$k_1$\emph{\ of }$l_1$\emph{\
and a component }$k_2$\emph{\ of }$l_2$\emph{, each of which has index 0 or
2. The index of the composed component }$k_1\#k_2$\emph{\ is equal to either
}$i(k_1) $\emph{\ or }$i(k_2)$\emph{. Finally, }$m$\emph{\ is a meridian of }%
$k_1\#k_2$\emph{\ with }$i=1.$

5) $V(l_1):$\emph{\ Choose a component }$k_1$\emph{\ of }$l_1$\emph{\ of
index 0 or 2, and replace N(k}$_1$\emph{,$S$}$^3$\emph{) by }$D^2\times S^1$%
\emph{\ with three indexed circles in it; \{0\}}$\times $\emph{\ }$S^1,$%
\emph{\ }$k_2$\emph{\ and }$k_3$\emph{. Here, }$k_2$\emph{\ and }$k_3$\emph{%
\ are parallel $(p,q)$-cables on }$\partial N($\emph{\{0\}}$\times $\emph{\ }%
$S^1,$\emph{\ }$D^2\times S^1)$\emph{, where p is the number of longitudinal
turns and q the number of the transverse ones}$.$\emph{\ The indices of \{0\}%
}$\times $\emph{\ }$S^1$\emph{\ and }$k_2$\emph{\ are either 0 or 2, and one
of them is equal to }$i(k_1)$\emph{. The index of }$k_3$\emph{\ is 1.}

6) $VI(l_{1})$\emph{: Choose a component }$k_{1}$\emph{\ of }$l_{1}$\emph{\
of index 0 or 2. Replace N(k}$_{1}$\emph{, $S$}$^{3}$\emph{) by }$%
D^{2}\times S^{1}$\emph{\ with two indexed circles in it; \{0\}}$\times $%
\emph{\ }$S^{1}$\emph{\ and the }$\emph{(2,q)}$\emph{-cable }$k_{2}$\emph{\
of \{0\}}$\times $\emph{\ }$S^{1}$\emph{. The index of \{0\}}$\times $\emph{%
\ }$S^{1}$\emph{\ is 1, and }$\emph{i}$\emph{(}$k_{2}$\emph{)=}$\emph{i}$%
\emph{(}$k_{1}$\emph{).}

Operations $I,$ $II$ and $III$ defined by Wada involve at least one
inessential circle of attachment; so, the flows characterized only by these
three operations have all their saddle orbits unknotted and unlinked. We
refer to them as type $A$ operations and denote by $\mathcal{F}_{A}\left(
S^{3}\right) $ the set of such flows and by $\mathcal{L}_{A}\left(
S^{3}\right) $ the set of the corresponding links of periodic orbits.

Operations $IV,$ $V$ and $VI$ imply only essential attachments; so, saddle
orbits generated by these operations are linked to other periodic orbits and
they can be knotted.

We call \emph{basic flows} those obtained from one operation of Wada on Hopf
links.

Let us observe that even though the round handle decomposition of a flow
is unique except for commutativity of some of the attachments involved,
non-equivalent flows can be characterized by the same link.

In this paper, we focus\ on the set of flows $\mathcal{F}_{A}\left(
S^{3}\right) $. We reproduce these flows by identifying two fat round
handles along their boundaries and we obtain that there exists an order in
the saddle orbits when heteroclinic trajectories appear.

\section{\bf Fat handles for $\mathcal{F}_{A}$ flows on S$^{3}$ \label{seccion asas}}

The 3-sphere $S^{3}$\ is composed of two solid tori identified along their
boundaries. So, we can obtain NMS flows on $S^{3}$ by identifying properly one repulsive and
one attractive tori along their boundaries. These tori with a flow pointing
inwards or outwards correspond to the fat round handles.

Given a flow $\varphi $ on $S^{3}$, a repulsive fat handle is obtained by
removing one attractive orbit and an attractive fat handle is obtained by
removing one repulsive orbit. In this section, we obtain and classify the
fat handles for $\mathcal{F}_{A}(S^{3})$ flows.

In the following we show the fat handles obtained by removing one orbit in
the basic flows of $\mathcal{F}_{A}(S^{3})$. We refer to them as \emph{basic
fat round handles} and we denote them by describing the periodic orbits that
contains. If the fat handle has an attractive or repulsive orbit in its core
we refer to it as thick torus; if it has no orbit in his core, we refer to it as solid torus.

Let $h$ denote the Hopf link, let $d$ denote a trivial separated
knot corresponding to an attractive or repulsive orbit, let $u$ denote a
trivial separated knot corresponding to a saddle orbit and let $\cdot$ denote
the separated sum of links.

\begin{figure}
\begin{center}
  \includegraphics[width=8cm]{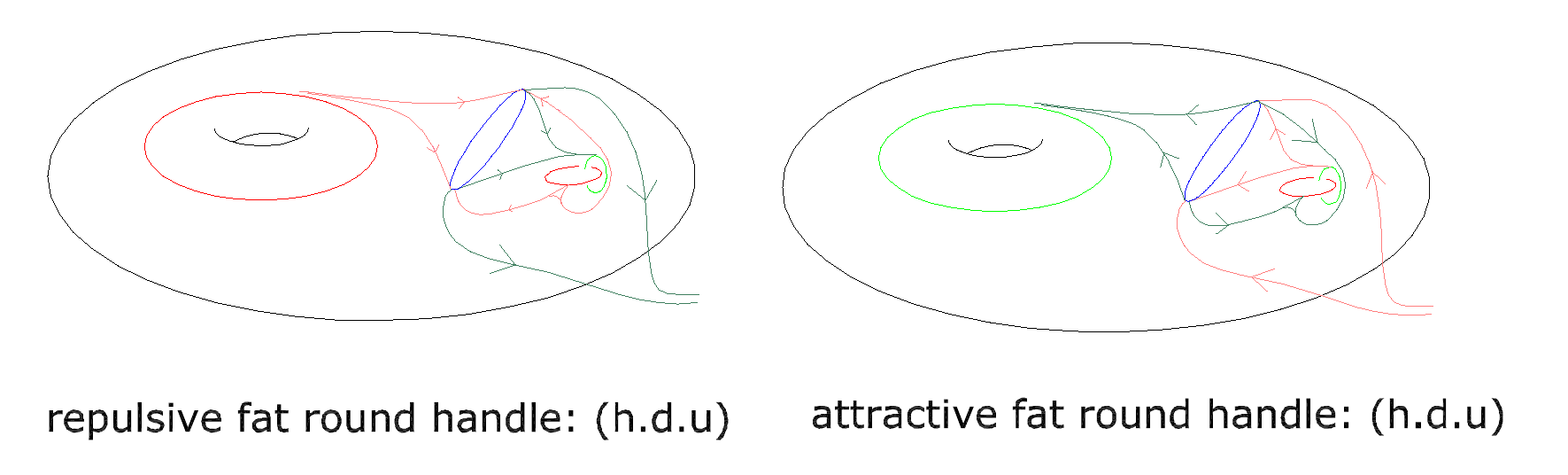}\\
  \caption{Fat handles associated to operation $I$  of Wada. }\label{hdu}
  \end{center}
\end{figure}

Consider the flow $I(h,h).$ The link of periodic orbits of this flow
consists in the separated sum of two Hopf links $h$ and an unknot $u$, $%
h\cdot h\cdot u$. We obtain the fat handles associated to operation $I$ by
removing one attractive or one repulsive orbit. These fat handles are tori
with an orbit in the core, i.e., thick tori, and we denote them $hdu$ (see Figure \ref{hdu}).

\begin{figure}[h]
\begin{center}
  \includegraphics[width=8cm]{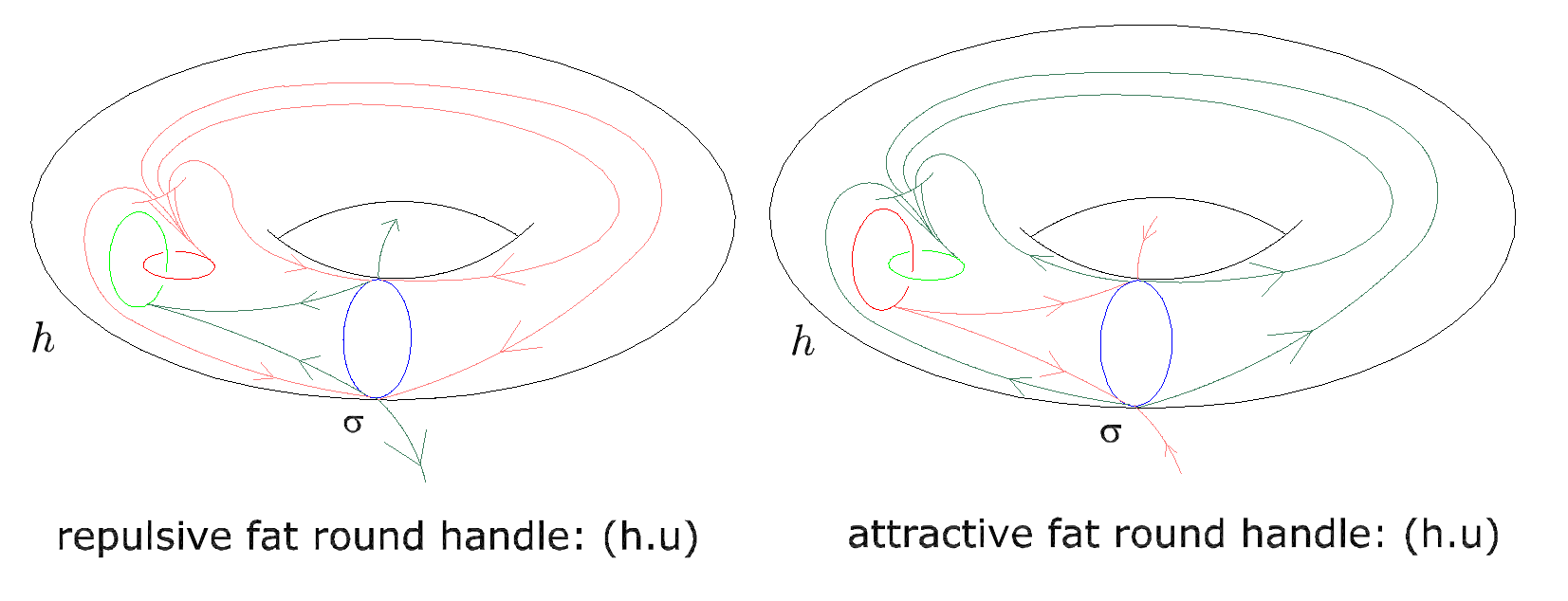}\\
  \caption{Fat handles associated to operation $II$ of Wada removing an orbit $d.$}\label{hu}
  \end{center}
\end{figure}

\begin{figure}[h]
\begin{center}
  \includegraphics[width=8cm]{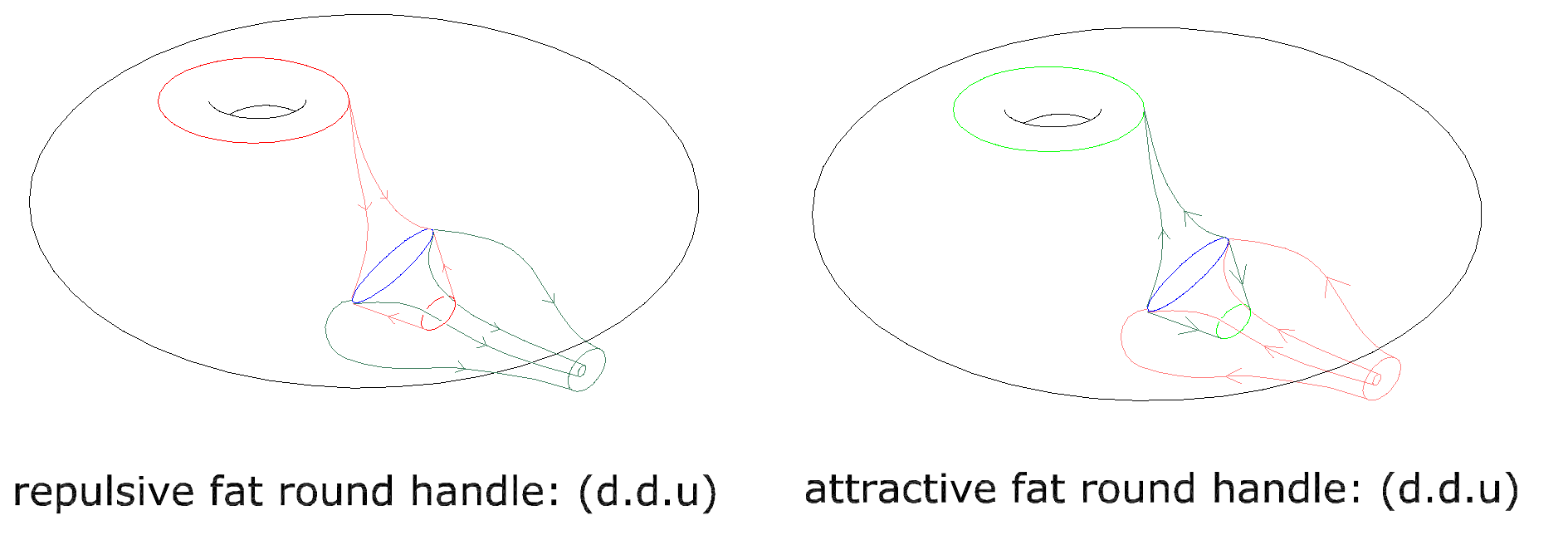}\\
  \caption{Fat handles associated to operation $II$ of Wada removing  a component of  $h.$}\label{ddu}
  \end{center}
\end{figure}

From the flow $II(h,h)=h\cdot d\cdot u$, we obtain the fat handles
associated to operation $II$ by removing one attractive or repulsive orbit.
Depending on the removed orbit, we obtain different types of fat handles. If
we take off the separated orbit $d,$ the fat handle is a torus without any
orbit in the core, i.e. solid torus (see Figure \ref{hu}) and we denote it
$hu$. If one component of the Hopf link $h$ is removed the fat handle is a
torus with one orbit in its core (see Figure \ref{ddu}) and we denote it $%
ddu $. In this last case, the orbits of type $d$\ can be both repulsive,
both attractive or one attractive and one repulsive.

\begin{figure}[h]
\begin{center}
  \includegraphics[width=8cm]{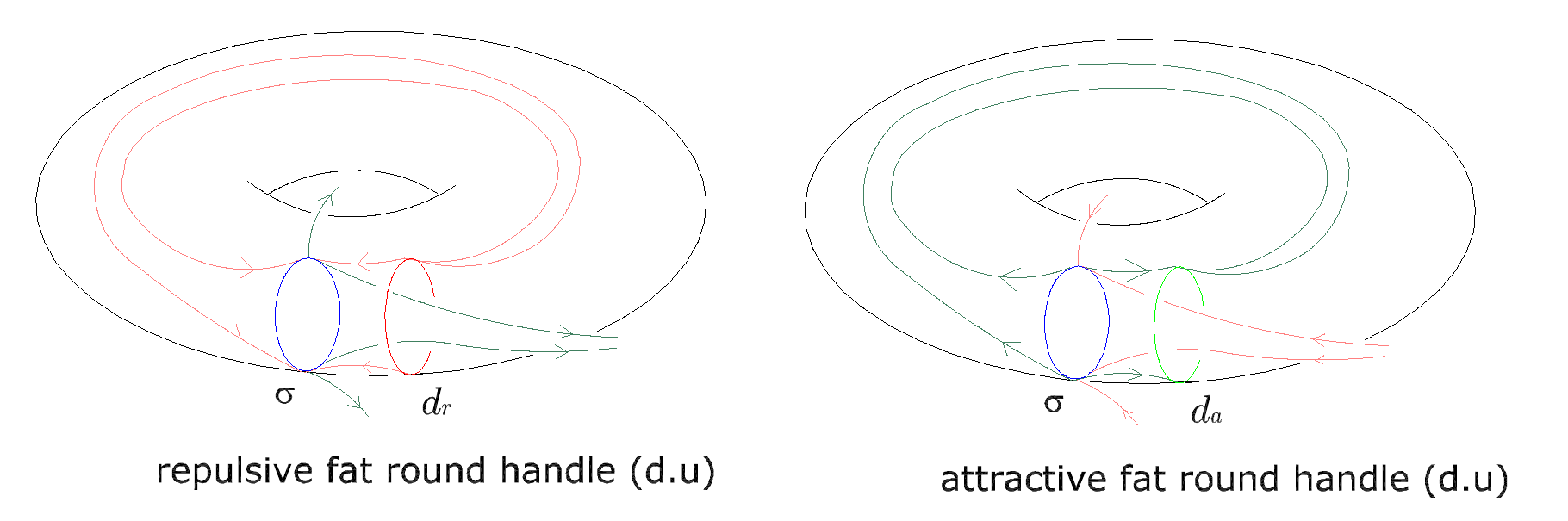}\\
  \caption{Fat handles associated to operation $III$ of Wada.}\label{du}
  \end{center}
\end{figure}

From the flow $III(h,h)=d\cdot d\cdot u$, we obtain the fat handles
associated to operation $III$ by removing one attractive or one repulsive
orbit. They are tori without any orbit in the core and we denote them $du$
(see Figure \ref{du}).

Let us remark that the identification along their boundaries of two fat
handles without any orbit in their cores yields a transversal intersection
of two invariant manifolds of saddle orbits (see \cite{campos2}).

We distinguish different classes of fat handles depending on the way the
invariant manifolds of the saddle leave or enter the torus and the number of
orbits in the canonical region where the identification of fat handles will
be made:

\begin{itemize}
\item The fat round handles that are thick tori and the invariant manifolds
of the saddle orbit going outwards the torus correspond to inessential
circles. The attractive (repulsive) orbits filling the essential hole of the
thick tori are in the canonical region where the identification with another
fat handle will be done (see Figures \ref{hdu} and \ref{ddu}).

\item The fat round handles that are solid tori coming from operation $II$:
one invariant manifold of the saddle orbit going outwards the torus is an
essential circle and there is not any attractive or repulsive orbit in the
canonical region of the identification (see Figure \ref{hu}).

\item The fat round handles that are solid tori coming from operation $III$:
one invariant manifold of the saddle orbit leaves or enters the torus
essentially and another one does it inessentially. There is one attractive
or repulsive orbit filling a non essential hole in the torus, in the
canonical region of the identification (see Figure \ref{du}).
\end{itemize}

Then, we define three classes of fat handles:

\begin{itemize}
\item A repulsive (attractive) fat handle belongs to class $\left[ I\right] $
if it corresponds to a thick torus and the invariant manifolds of the
saddles orbits go outwards (inwards) the torus by means of inessential
circles.

\item A repulsive (attractive) fat handle belongs to class $\left[ II\right]
$ if it corresponds to a solid torus, the invariant manifolds of the saddles
orbits go outwards (inwards) the torus by means of essential circles and
there is not any attractive or repulsive orbit in the canonical region of
the identification.

\item A repulsive (attractive) fat handle belongs to class $\left[ III\right]
$ if it corresponds to a solid torus, the invariant manifolds of the saddles
orbits go outwards (inwards) the torus by means of essential and inessential
circles and there is one attractive or repulsive orbit, filling a non
essential hole in the torus, in the canonical region of the identification.
\end{itemize}

From the identification along their boundaries of one attractive and one
repulsive basic fat handles we obtain the $\mathcal{F}_{A}(S^{3})$ flows
with two saddle orbits. Iterated fat round handles with two saddles are
obtained by removing one repulsive (attractive) orbit in these flows.

Let us notice that when two fat handles belonging to class $\left[ I\right] $
are identified along their boundaries, the orbits in their cores form a Hopf
link in the canonical region of the identification. On the other hand, when
fat handles of class $\left[ II\right] $ and $[III]$ are identified, one
heteroclinic orbit connecting the two saddles orbits appears (see Figures %
\ref{22}, \ref{23} and \ref{33}).

In the following proposition we show that the iterated fat handles with two
saddle orbits can be classified in one of these three classes of fat handles
defined above.

\begin{proposition}
\label{2us} For $ \mathcal{F}_{A}\left( S^{3}\right)-$flows, a fat handle with two saddle orbits belongs to class $\left[ I
\right] ,$ $\left[ II\right] $ or $\left[ III\right] .$
\end{proposition}

\noindent
\proof The fat handles with two saddles are obtained by identifying one
repulsive and one attractive basic fat handles along their boundaries and
then removing a repulsive or attractive periodic orbit.

\begin{itemize}
\item When an attractive fat handle $hdu$ is identified with a repulsive fat
handle $hdu,$ the flow $I\left( I\left( h,h\right) ,h\right) =h\cdot h\cdot
h\cdot u\cdot u$ is obtained$.$ If an attractive (repulsive) orbit is
removed after the identification, the resulting repulsive (attractive) fat
handle, $hhduu$, corresponds to a thick torus with the manifolds of the
saddle orbits going outwards (inwards) the torus inessentially. So, it
belongs to class $\left[ I\right] .$

\item If an attractive fat handle $hdu$ is identified with a repulsive fat
handle $ddu,$ the flow $I\left( II\left( h,h\right) ,h\right) =h\cdot h\cdot
d\cdot u\cdot u$ is obtained. Now, we can remove the orbit $d$ or a
component of a Hopf $h.$ If a component of $h$ is removed, the repulsive
(attractive) fat handle is $hdduu$,\ with one toroidal and one Hopf holes
inside and invariant manifolds of the saddles crossing unessentially the
torus; so, it belongs to class $\left[ I\right] $. If the orbit $d$ is
removed, the fat handle is $hhuu,$ a solid torus with two Hopf holes inside,
with no orbit in the canonical region of the identification and with the
invariant manifold of one of the saddle orbits crossing the torus by means
of an essential circle; so, it belongs to class $\left[ II\right] $.

\item The identification of an attractive (repulsive) fat handle $hdu$ with
a repulsive (attractive) fat handle $hu$ is not admissible because this
identification generates a bitorus (see Figure \ref{bitoro}).

\begin{figure}[h]
\begin{center}
  \includegraphics[width=6cm]{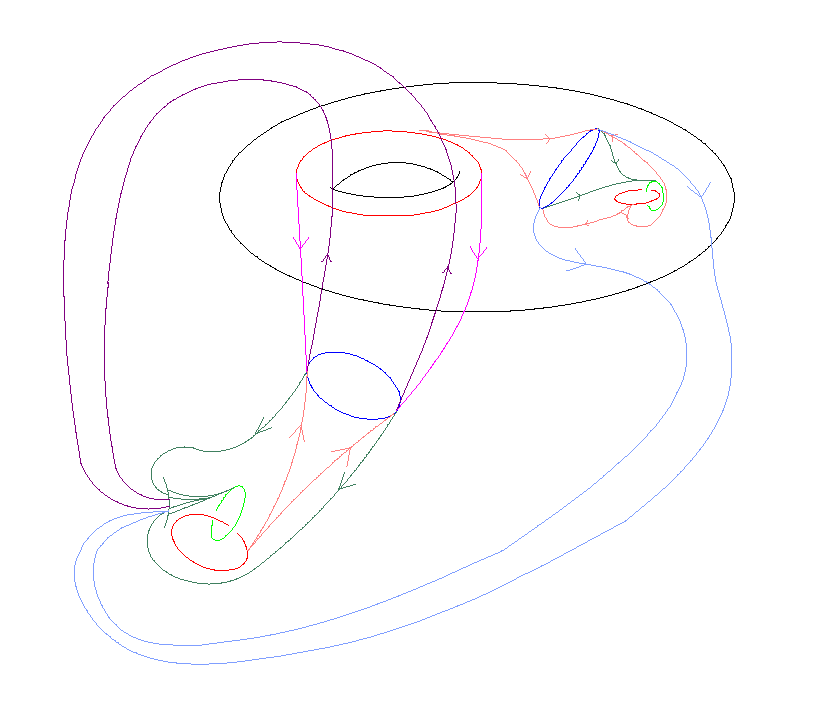}\\
  \caption{Identification of fat handles $hdu$ and $hu$ along their boundaries.}\label{bitoro}
  \end{center}
\end{figure}

\item When an attractive fat handle $hdu$ is identified with a repulsive fat
handle $du,$ we obtain the flow $I\left( III\left( h,h\right) ,h\right)
=h\cdot d\cdot d\cdot u\cdot u.$ If an attractive (repulsive) orbit is now
removed, the resulting repulsive (attractive) fat handle is either a thick
torus with two toroidal holes inside, $ddduu,$ belonging to class $\left[ I%
\right] $ or a solid torus $hduu$ with one orbit in the canonical region of
the identification, belonging to class $\left[ III\right] $.

\item If an attractive fat handle $ddu$ is identified with a repulsive fat
handle $ddu,$ the resulting flow is $II\left( II\left( h,h\right) ,h\right)
=h\cdot d\cdot d\cdot u\cdot u.$ Now, we can remove a component of the Hopf $%
h$ or an orbit $d.$ Then, the resulting repulsive (attractive) fat handle is
either a thick torus with two toroidal holes inside, $ddduu,$ belonging to
class $\left[ I\right] $ or a solid torus with one toroidal and one Hopf
holes inside, $hduu,$ with no orbit in the canonical region of the
identification, belonging to class $\left[ II\right] .$

\item The identification of an attractive (repulsive) fat handle $ddu$ with
a repulsive (attractive) fat handle $hu$ is not admissible because a bitorus
is obtained.

\item When an attractive fat handle $ddu$ is identified with a repulsive fat
handle $du$, the defined flow is $II\left( III\left( h,h\right) ,h\right)
=d\cdot d\cdot d\cdot u\cdot u$. If an attractive (repulsive) orbit is
removed after the identification the resulting repulsive (attractive) fat
handle $dduu$ is either a solid torus with two toroidal holes inside and
with no orbit in the canonical region of the identification belonging to
class $\left[ II\right] $ or a solid torus with two toroidal holes inside
with one orbit $d$ in the canonical region of the identification, belonging
to class $\left[ III\right] .$

\item If an attractive fat handle $hu$ is identified with a repulsive fat
handle $hu,$ the flow$II\left( II\left( h,h\right) ,h\right) =$ $h\cdot
h\cdot u\cdot u$ is obtained (see Figure \ref{22}). In this case, one
heteroclinic trajectory appears. If a repulsive (attractive) orbit is
removed after the identification the resulting attractive (repulsive) fat
handle $dhuu$ is a thick torus and it belongs to class $\left[ I\right] $.%

\begin{figure}[h]
\begin{center}
  \includegraphics[width=6cm]{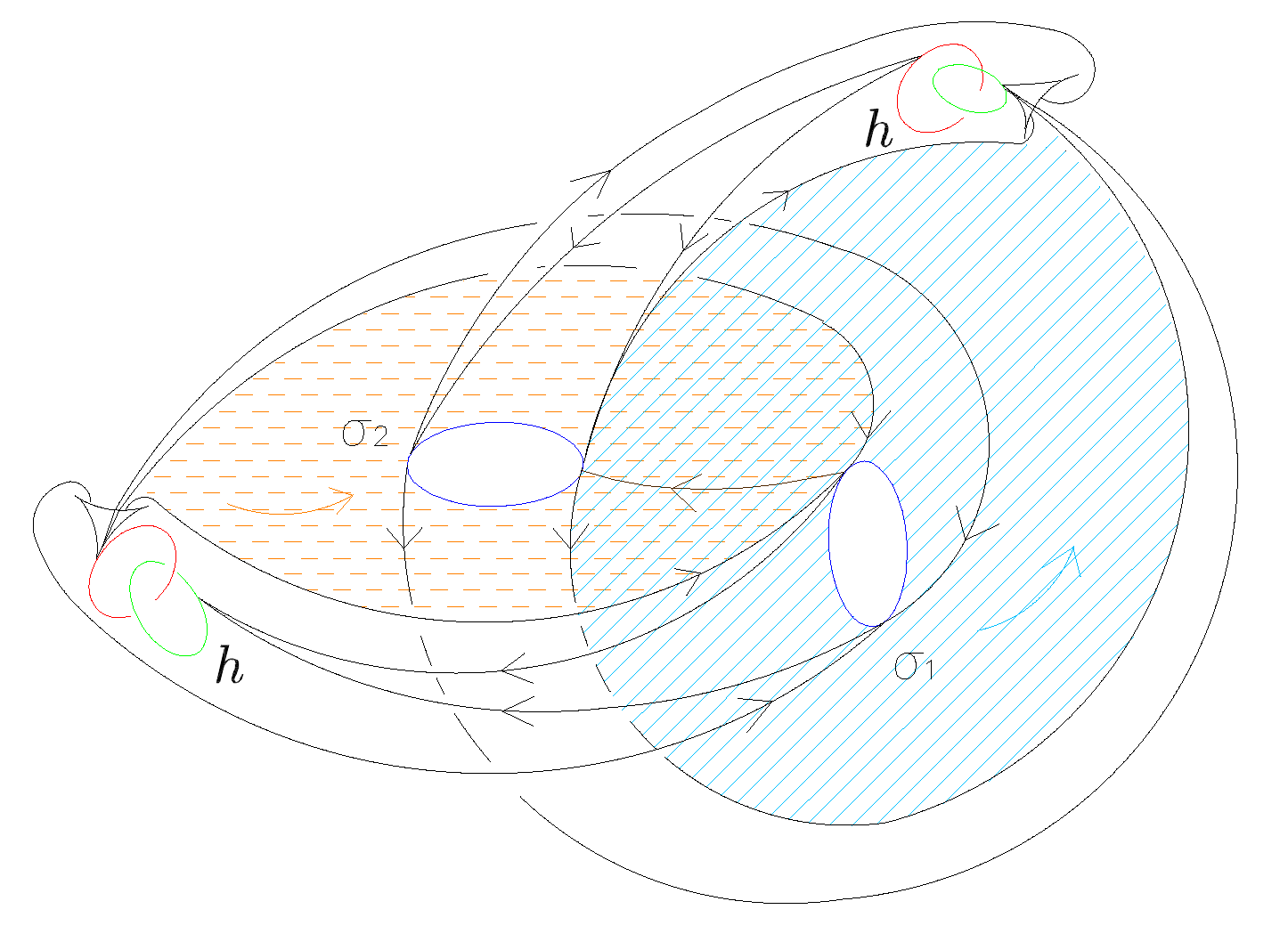}\\
  \caption{Flow $II\left( II \left(h,h\right) ,h\right) =h\cdot h\cdot u\cdot u.$}\label{22}
  \end{center}
\end{figure}

\item If an attractive (repulsive) fat handle $hu$ is identified with a
repulsive (attractive) fat handle $du,$ the resulting flow is $II\left(
III\left( h,h\right) ,h\right) =h\cdot d\cdot u\cdot u$ and one heteroclinic
trajectory appears (see Figure \ref{23}). If a repulsive (attractive) orbit
is removed after the identification the resulting attractive (repulsive) fat
handle is either a solid torus with one Hopf hole inside, $huu,$ belonging
to class $\left[ II\right] $, or a thick torus with one toroidal hole
inside, $dduu,$ belonging to class $\left[ I\right]$ .

\begin{figure}[h]
\begin{center}
  \includegraphics[width=6cm]{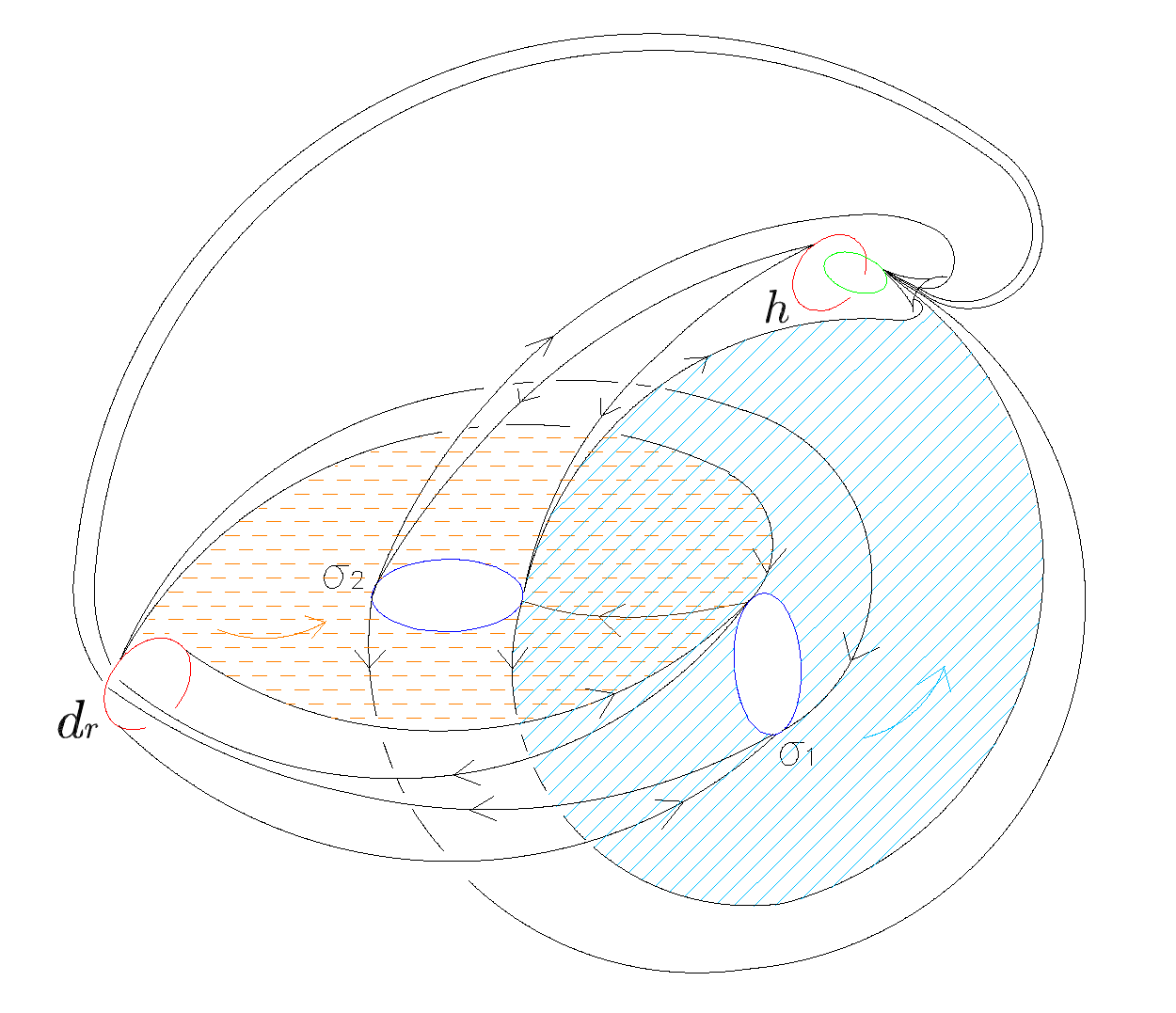}\\
  \caption{Flow $II\left( III\left( h,h\right) ,h\right) =h\cdot d\cdot u\cdot u$}\label{23}
  \end{center}
\end{figure}

\item Finally, if an attractive fat handle $du$ is identified with a
repulsive fat handle $du,$ we obtain the flow $III\left( III\left(
h,h\right) ,h\right) =d\cdot d\cdot u\cdot u$ (see Figure \ref{33}).
If an attractive (repulsive) orbit is removed after the identification
the resulting repulsive (attractive) fat handle $duu$ is a solid torus with
one orbit $d$ in the canonical region of the identification; so, it belongs
to class $\left[ III\right] $.

\begin{figure}[h]
\begin{center}
  \includegraphics[width=6cm]{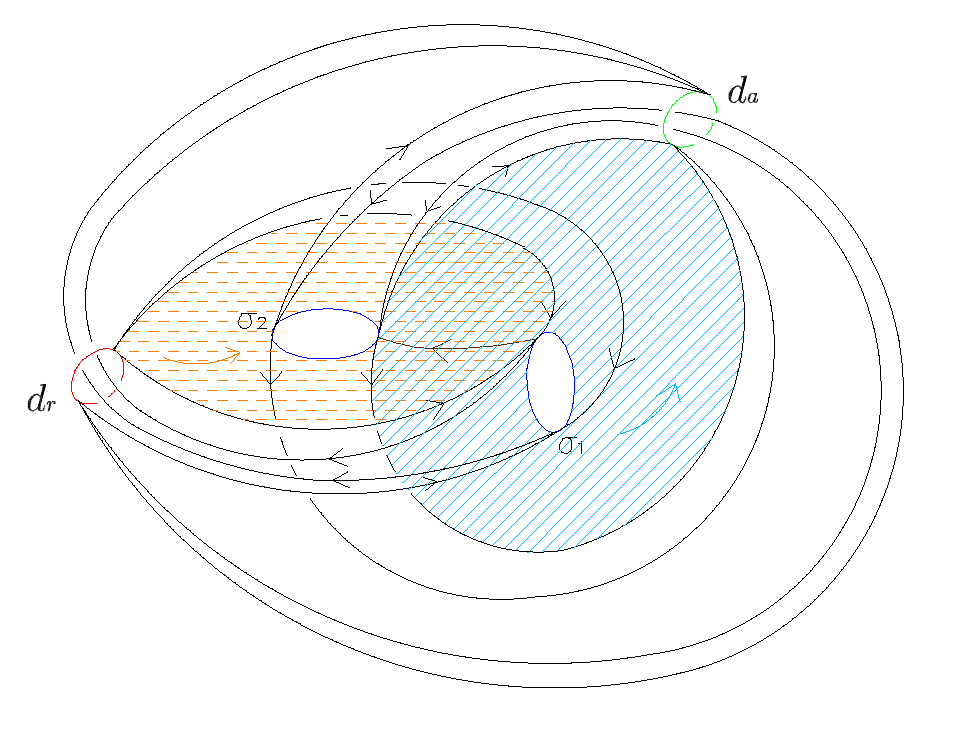}\\
  \caption{Identification of two fat handles $du.$}\label{33}
  \end{center}
\end{figure}

\end{itemize}

\endproof

Let us remark that fat handles of class $\left[ I\right] $ are of the form $%
ldu^{2},$where $d$ is in the core of the fat handle and $l$ represent the
other orbits. Similarly, fat handles of class $\left[ II\right] $ are of the
form $lhu^{2},$where $h$ is in the canonical region of the identification
and fat handles of class $\left[ III\right] $ are of the form $lu^{2},$where
$l$ contains only one $d$ in the region of identification.

Along the proof of\ the previous proposition we have obtained all the flows $%
\varphi $\ with two unlinked saddles. Therefore,

\begin{proposition}
A flow $\varphi \in \mathcal{F}_{A}\left( S^{3}\right) $\ with two saddle
orbits can be obtained by identifying one repulsive and one attractive basic
fat handles along their boundaries.
\end{proposition}

\proof From Proposition \ref{2us} we have that:

The flow $I\left( I\left( h,h\right) ,h\right) =h\cdot h\cdot h\cdot u\cdot
u $ is obtained by identifying the fat handles $hdu$ and $hdu\ $along their
boundaries.

The flow $I\left( II\left( h,h\right) ,h\right) =h\cdot h\cdot d\cdot u\cdot
u$ is obtained by identifying the fat handles $hdu$ and $ddu$ along their
boundaries.

The flow $I\left( III\left( h,h\right) ,h\right) =h\cdot d\cdot d\cdot
u\cdot u$ is obtained by identifying the fat handles $hdu$ and $du\ $along
their boundaries.

The flow $II\left( II\left( h,h\right) ,h\right) =h\cdot d\cdot d\cdot
u\cdot u$ is obtained by identifying the fat handles $ddu$ and $ddu\ $along
their boundaries.

The flow $II\left( II\left( h,h\right) ,h\right) =h\cdot h\cdot u\cdot u$ is
obtained by identifying the fat handles $hu$ and $hu\ $along their
boundaries.

The flow $II\left( III\left( h,h\right) ,h\right) =d\cdot d\cdot d\cdot
u\cdot u$ is obtained by identifying the fat handles $ddu$ and $du\ $along
their boundaries.

The flow $II\left( III\left( h,h\right) ,h\right) =h\cdot d\cdot u\cdot u$
is obtained by identifying the fat handles $hu$ and $du\ $along their
boundaries.

The flow $III\left( III\left( h,h\right) ,h\right) =d\cdot d\cdot u\cdot u$
is obtained by identifying the fat handles $du$ and $du\ $along their
boundaries.
\endproof

Let us recall that there are non-equivalent flows characterized by the same
link of periodic orbits.

As we see in the following propositions, the fat handles with any number of
unlinked saddle orbits are in the same classes previously defined and the
flows $\mathcal{F}_{A}\left( S^{3}\right) $ can be obtained by replacing
successively one attractive (repulsive) orbit by the corresponding basic fat
handle, taking into account the restrictions in the proof of Proposition \ref%
{2us}: fat handles belonging to class $\left[ I\right] $ can not be
identified with fat handles belonging to class $\left[ II\right] ,$ i.e.,
one component of a Hopf link can not be replaced by a fat handle belonging
to class $\left[ II\right] ${\Large .}

\begin{proposition}
\label{nus}  For $ \mathcal{F}_{A}\left( S^{3}\right)-$flows, a fat handle with $n$ saddle orbits belongs to class $\left[ I%
\right] ,$ $\left[ II\right] $ or $\left[ III\right] .$
\end{proposition}

\proof Let us show this result by means of an induction process, by
identifying iterated fat handles with the basic ones.

We know that the first identification leads to fat handles with $n=2$
saddles and they belong to class $[I],$ $[II]$ or $[III].$ Let us suppose
that it holds for fat handles with $n-1$ saddles and let us prove it for $n$.

Let $l\cdot u^{n-1}$ denote the link corresponding to a flow $\varphi $ with
$n-1$ saddle orbits. By assumption, after removing one attractive or
repulsive orbit $k,$ the corresponding fat handle is of class $\left[ I%
\right] ,$ $\left[ II\right] $ or $\left[ III\right] $, \ depending on the
removed orbit $k$. \ According to the previous notation, we denote them as $%
\left( l-h\right) du^{n-1}$, $\left( l-h-k\right) hu^{n-1}$ or $\left(
l-k\right) u^{n-1}$, respectively.

Now, we make the identifications with the different basic fat
handles and then, we remove an attractive or repulsive periodic orbit in
order to obtain the fat handles with $n$ saddles$.$ As we identify repulsive
with attractive fat handles, the removed orbit is located in the canonical
region where\ the identification is made; the other canonical regions do not
change.

\begin{itemize}
\item Let us suppose that the repulsive fat handle with $\left( n-1\right) $
saddle orbits belongs to class $\left[ I\right] .$ It means that it is a
thick torus with a repulsive orbit filling the essential hole and all the
manifolds of the saddles go outwards the torus in an inessential way. Let us
denote it $\left( l-h\right) du^{n-1},$ where $d$ fills the core of the
thick torus. If this fat handle is identified with an attractive fat handle $%
hdu,$ the attractive orbit $d$ of this fat handle forms a Hopf link with the
orbit in the core of the repulsive fat handle and the flow $l_{2}=l\cdot
h\cdot u^{n}=l_{1}\cdot h\cdot u=I(l_{1},h)$ is obtained.

If a repulsive (attractive) orbit belonging to a Hopf link is removed after
the identification the resulting attractive (repulsive) fat handle is also a
thick torus and it belongs to class $\left[ I\right] .$

\item Let us suppose that the attractive fat handle $\left( l-h\right)
du^{n-1}$ of class $\left[ I\right] $ is identified with a repulsive fat
handle $ddu$ of class $\left[ I\right] .$
The two orbits that are in the core of the thick tori become a Hopf link
after identifying them along their boundaries and the resulting flow is $%
l_{2}=l\cdot d\cdot u^{n}=II\left( l_{1},h\right) ,$ where $d$ is an
attractive or repulsive orbit filling a non-essential toroidal hole. If it
is identified with the basic fat handle $du$ (class $\left[ III\right] ),$
the resulting flow is $l_{2}=\left( l-h\right) \cdot d\cdot d\cdot
u^{n}=\left( l-k\right) \cdot d\cdot u^{n},$ where each $d$ is an attractive
or repulsive orbit filling a non-essential toroidal hole. As in the previous
proposition, depending on the position of the removed orbit, the resulting
fat handle is in class $\left[ I\right] ,$ $\left[ II\right] $ or $\left[ III%
\right] .$

\item An attractive fat handle of class $\left[ I\right] $ can not be
identified with a repulsive fat handle $hu$ because a bitorus appears.

\item An attractive fat handle of type $\left[ II\right] ,$ $\left(
l-h-k\right) hu^{n-1}$ can be identified with a repulsive fat handle $hu.$
The flow corresponds to $\left( l-k\right) \cdot h\cdot u^{n}$ and one
heteroclinic trajectory appears in the canonical region where the
identification has been made; this canonical region becomes a 3-ball. After
removing one orbit belonging to a Hopf link, the resulting fat handle is in
class $\left[ I\right] .$

\item If an attractive fat handle of type $\left[ II\right] ,$ $\left(
l-h-k\right) hu^{n-1}\ $is identified with a repulsive fat handle $du,$ the
resulting flow corresponds to $\left( l-k\right) \cdot d\cdot u^{n}$ and one
heteroclinic trajectory appears. If the repulsive (attractive) orbit removed
is in the canonical region where the identification has been made, the fat
handle is a solid torus belonging to class $\left[ II\right] $.

\item Finally, if an attractive fat handle of type $\left[ III\right] ,$ $%
\left( l-k\right) u^{n-1},$ is identified with a repulsive fat handle $du,$
the resulting flow corresponds to $\left( l-k\right) \cdot d\cdot u^{n}$ and
one heteroclinic trajectory appears. If the repulsive (attractive) orbit
removed is in the canonical region where the identification has been made,%
{\large \ }the fat handle is a solid torus belonging to class $\left[ III%
\right] $. \endproof
\end{itemize}

Let us emphasize that, along the proof, we have obtained all the flows $%
\varphi $\ with $n$ unlinked saddles. Therefore,

\begin{theorem}
\label{lnus} A flow $\varphi \in \mathcal{F}_{A}\left( S^{3}\right) $\ with $%
n$ saddle orbits can be obtained by identifying fat handles along their
boundaries.
\end{theorem}

\proof Following the results of Proposition \ref{nus}, if $l\cdot u^{n-1}$
is a flow with $n-1$ saddle orbits, the corresponding fat handle obtained by
removing one attractive (or repulsive) orbit is one of the defined in the
previous proposition: $\left( l-h\right) du^{n-1},$ $\left( l-h-k\right)
hu^{n-1}$ and $\left( l-k\right) u^{n-1}$ if it is of class $\left[ I\right]
,$ $\left[ II\right] $ or $\left[ III\right] ,$ respectively$.$

Iterating one more Wada operation, we have: $I\left( lu^{n-1},h\right)
=l\cdot h\cdot u^{n},$ $II\left( lu^{n-1},h\right) =l\cdot d\cdot u^{n}$ or $%
II\left( lu^{n-1},h\right) =\left( l-k\right) \cdot h\cdot u^{n}$ and $%
III\left( lu^{n-1},h\right) =\left( l-k\right) \cdot d\cdot u^{n},$ that are
the flows obtained in the Proposition \ref{nus}. \endproof

Let us notice that one heteroclinic trajectory appears when fat handles that
are solid tori are identified along their boundaries. As we see in the
following section, these type of trajectories impose an order in the flow.

\section{\bf Order in the NMS flows \label{sectorder}}

The round handle decomposition of the manifold $S^{3}$ for a given NMS flow
is unique except for commutativity of some of the attachments involved. There
exist attachments that can lead to different flows depending on the order
they are made. The appearance of heteroclinic trajectories connecting saddle
orbits implies non commutativity of the operations involved establishing an
order in the round handle decomposition. As we have seen in the previous
section, heteroclinic trajectories connecting saddles appear with the
identification of fat handles of classes $\left[ II\right] $ and $\left[ III%
\right] .$

\begin{proposition}
\label{orden} For $\varphi \in \mathcal{F}_{A}(S^{3})$, the heteroclinic
trajectories induce an order in the flow$.$
\end{proposition}

\begin{figure}[h!]
\begin{center}
  \includegraphics[width=12cm]{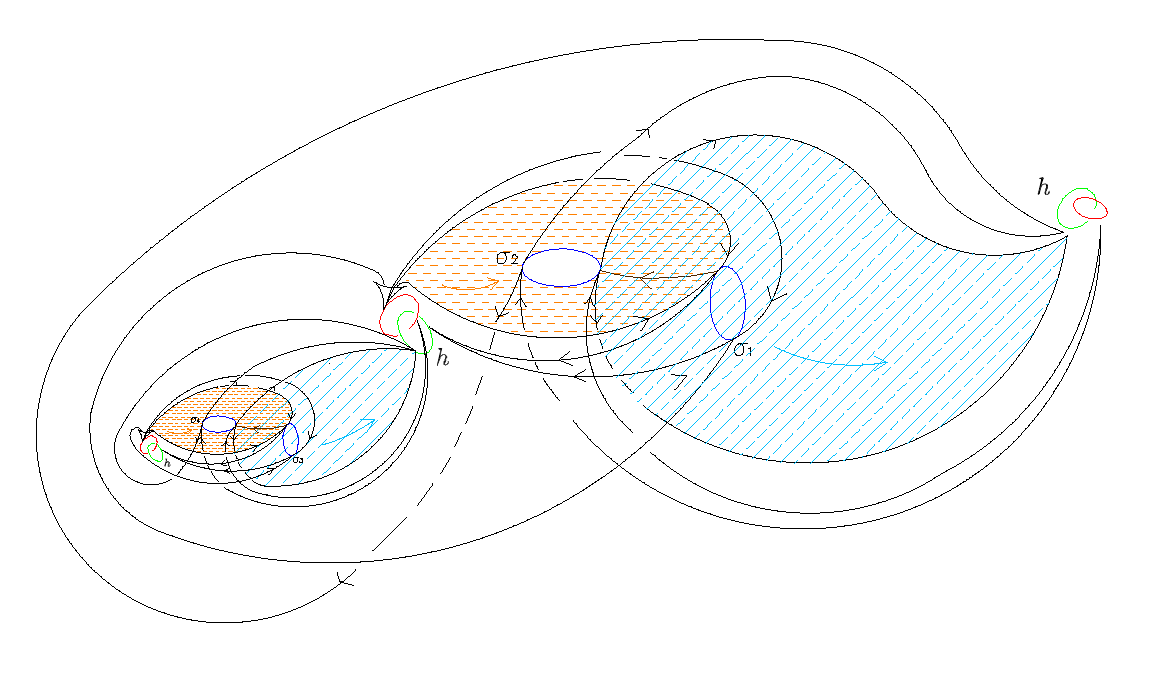}\\
  \caption{Flow $II(II(II(II(h,h),h),h),h)=h \cdot h\cdot h\cdot u\cdot u\cdot u\cdot u.$ There is a partial order $\sigma _{1}< \sigma _{2}$, $\sigma _{3}<\sigma _{4}.$} \label{2222}
  \end{center}
\end{figure}

\begin{figure}[h!]
\begin{center}
  \includegraphics[width=8cm]{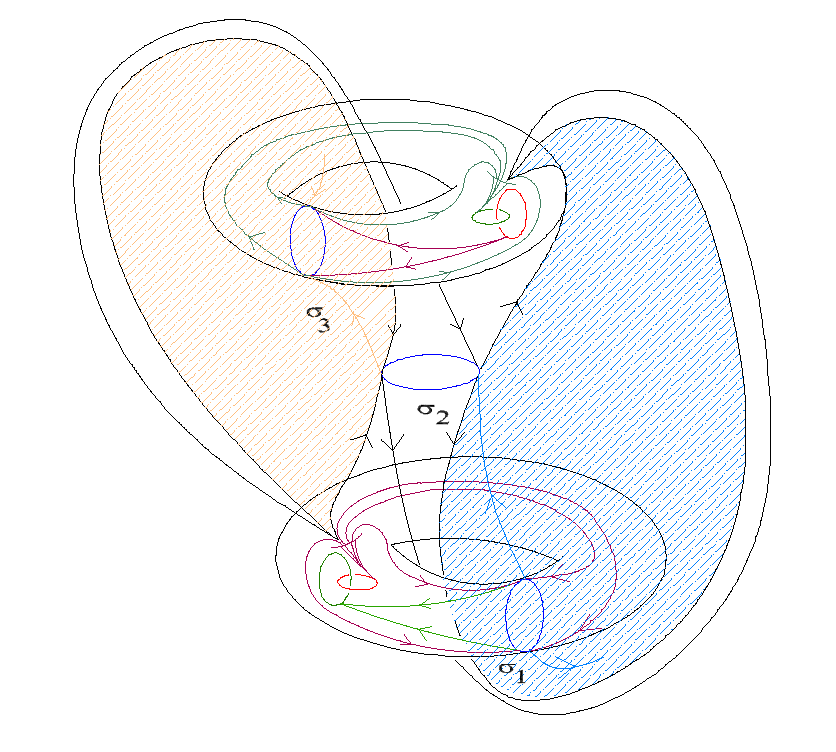}\\
  \caption{Flow $II(II(III(h,h),h),h)=h\cdot h\cdot u\cdot u\cdot u$. There is an order $\sigma _{1}<\sigma _{2}<\sigma _{3}.$}\label{223}
  \end{center}
\end{figure}

\proof Heteroclinic trajectories appear when fat round handles belonging to class $%
\left[ II\right] $ or $\left[ III\right] $ are identified along their
boundaries. The unstable manifold of a saddle orbit $u_{i}$ in the repulsive
fat handle intersects transversely the stable manifold of the saddle orbit $%
u_{i+1}$ in the attractive fat handle (see Figures \ref{2222} and \ref{223}).

Then, the saddles are ordered, $u_{i}<u_{i+1},$ and this order of the orbits
implies a partial order in the filtration of the flow.

\endproof

For the particular case of a flow where only operation $III$ is implied the
order is total.

\begin{theorem}
\label{orden3} Let $\mathcal{F}_{3}(S^{3})\ $be the set of NMS flows on $%
S^{3}$ coming only from operation $III.$ Then, the set of orbits of $\varphi
\in \mathcal{F}_{3}(S^{3})$\ is totally ordered.
\end{theorem}

\proof A flow $\varphi \in \mathcal{F}_{3}(S^{3})$ is of the form $d_{a}\cdot
d_{r}\cdot u\cdot ...\cdot $ $u,$ where $u$ denotes an unknot corresponding
to a saddle orbit. The flow goes from the repulsive orbit $d_{r}$ to the
attractive orbit $d_{a}.$

\begin{figure}[h!]
\begin{center}
  \includegraphics[width=8cm]{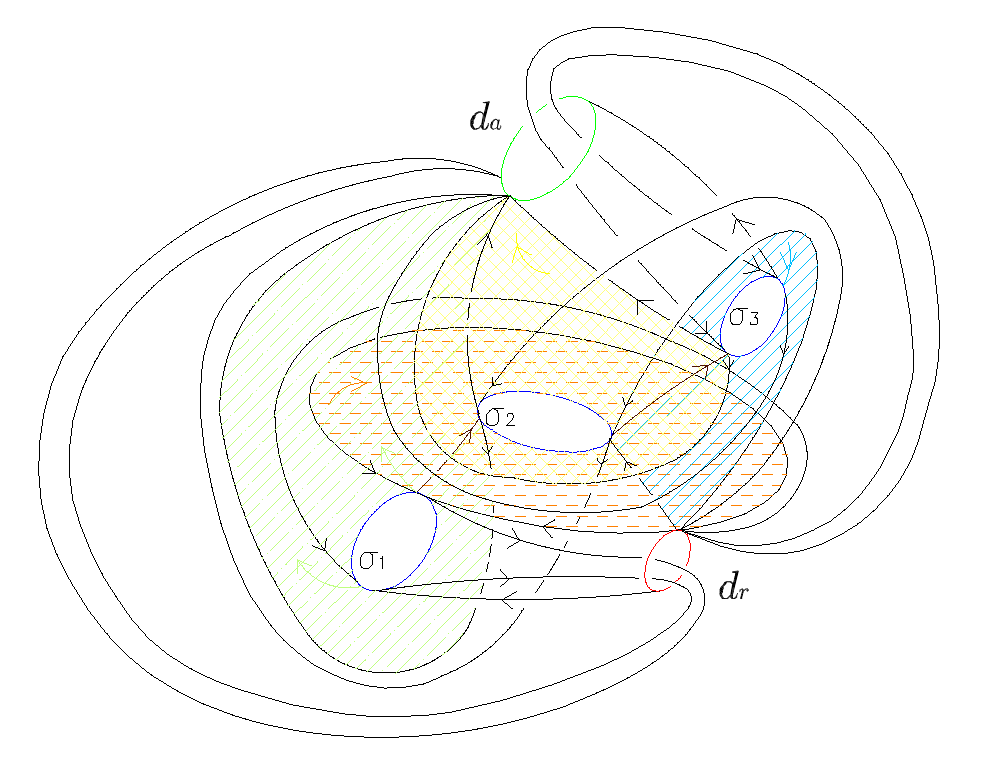}\\
  \caption{Flow $III(III(III(h,h),h),h).$ The order is $d_{r}<u_{1}<u_{2}<u_{3}<d_{a}.$}\label{333}
  \end{center}
\end{figure}

Consider the flow $III\left( III\left( h,h\right) ,h\right) =d_{a}\cdot
d_{r}\cdot u\cdot u$ (see Figure \ref{33}). One unstable manifold of the
first saddle orbit $u_{1}$ cuts transversely the stable manifold of the
second saddle orbit $u_{2}$ and one heteroclinic trajectory appears from $%
u_{1}$ to $u_{2}$ (see \cite{campos2}). We write it as $u_{1}<u_{2}.$

Each time operation $III$ is applied, two fat handles that are solid tori
are identified along their boundaries; so, the unstable manifold of the
saddle $u_{i}$ cuts transversely the stable manifold of the new saddle $%
u_{i+1}$ and one heteroclinic trajectory appears connecting $u_{i}$ and $%
u_{i+1}$; so, $u_{i}<u_{i+1}$ (see Figure \ref{333}).

Therefore, for a flow $d_{a}\cdot d_{r}\cdot u\cdot ...\cdot $ $u$ we can
write
\begin{equation}
d_{r}<u_{1}<...<u_{i}<u_{i+1}<...<u_{n}<d_{a}.\hspace{2cm}%
\end{equation}

\endproof

\end{document}